\documentclass{amsart}
\usepackage{amssymb}

\begin{document}

\def\uL(#1){\L^{\le #1}}

\def\CC{\mathbb C}
\def\NN{\mathbb N}
\def\RR{\mathbb R}
\def\TT{\mathbb T}
\def\ZZ{\mathbb Z}

\def\Bb{{\mathcal B}}
\def\Ee{{\mathcal E}}
\def\Ff{{\mathcal F}}
\def\Hh{{\mathcal H}}
\def\Kk{{\mathcal K}}
\def\Mm{{\mathcal M}}
\def\Oo{{\mathcal O}}
\def\Ss{{\mathcal S}}
\def\Tt{{\mathcal T}}

\def\Ad{\operatorname{Ad}}
\def\Aut{\operatorname{Aut}}
\def\clsp{\overline{\lsp}}
\def\End{\operatorname{End}}
\def\id{\operatorname{id}}
\def\Isom{\operatorname{Isom}}
\def\Ker{\operatorname{Ker}}
\def\lsp{\operatorname{span}}
\def\fin{\operatorname{fin}}
\def\Hom{\operatorname{Hom}}
\def\Obj{\operatorname{Obj}}
\def\dom{\operatorname{dom}}
\def\cod{\operatorname{cod}}
\def\FE{\operatorname{FE}}
\def\Ext{\operatorname{Ext}}
\def\MCE{\operatorname{MCE}}
\def\Mor{\operatorname{Mor}}

\def\Lmin{{%
\L^{\min}}}%

\def\L{\Lambda}
\def\margincomment(#1){{%
\begin{flalign*}
\hskip-1.5cm \text{\Large\sc {#1}}&&&&
\end{flalign*} }}%

\theoremstyle{plain}
\newtheorem{theorem}{Theorem}[section]
\newtheorem*{theorem*}{Theorem}
\newtheorem*{prop*}{Proposition}
\newtheorem{cor}[theorem]{Corollary}
\newtheorem{lemma}[theorem]{Lemma}
\newtheorem{prop}[theorem]{Proposition}
\newtheorem{conj}[theorem]{Conjecture}
\theoremstyle{remark}
\newtheorem{rmk}[theorem]{Remark}
\newtheorem{rmks}[theorem]{Remarks}
\newtheorem*{aside}{Aside}
\newtheorem*{note}{Note}
\newtheorem{comment}
[theorem]{Comment}
\newtheorem{example}[theorem]{Example}
\newtheorem{examples}[theorem]{Examples}
\theoremstyle{definition}
\newtheorem{dfn}[theorem]{Definition}
\newtheorem{dfns}[theorem]{Definitions}
\newtheorem{notation}[theorem]{Notation}

\numberwithin{equation}{section}

\title[Relative Cuntz-Krieger algebras]{Relative Cuntz-Krieger algebras of finitely aligned higher-rank graphs}
\author{Aidan Sims}
\address{Department of Mathematics  \\
      University of Newcastle\\
NSW  2308\\
AUSTRALIA}
\email{aidan@frey.newcastle.edu.au}
\date{December 8, 2003}
\thanks{This research is part of the author's PhD thesis, supervised by Professor Iain Raeburn, and was supported by an Australian Postgraduate Award and by the Australian Research Council.}

\begin{abstract} We define the relative Cuntz-Krieger algebras associated to finitely aligned higher-rank graphs. We prove versions of the gauge-invariant uniqueness theorem and the Cuntz-Krieger uniqueness theorem for relative Cuntz-Krieger algebras.
\end{abstract}

\maketitle

\section{Introduction}\label{sec:intro} 
Cuntz-Krieger algebras associated to directed graphs and their analogues have been of significant interest recently, due in large part to the explicit relationship between the loop-structure of a graph and the ideal structure of its Cuntz-Krieger algebra.

A directed graph $E$ consists of a collection $E^0$ of vertices, a collection $E^1$ of edges joining the vertices, and maps $r,s : E^1 \to E^0$ which indicate the ranges and sources of the edges. The Cuntz-Krieger algebra of $E$, denoted $C^*(E)$, is the universal algebra generated by mutually orthogonal projections $\{p_v : v \in E^0\}$ and by partial isometries $\{s_e : e \in E^1\}$ with mutually orthogonal range projections such that for $e \in E^1$, we have $s^*_e s_e = p_{s(e)}$ and such that whenever $v \in E^0$ satisfies $0 < |r^{-1}(v)| < \infty$, we have
\[
s_v = \sum_{e \in r^{-1}(v)} s_e s^*_e.
\]

The universal property of $C^*(E)$ ensures that it carries a strongly continuous \emph{gauge action} $\gamma$ of $\TT^k$ satisfying $\gamma_z(p_v) = p_v$ and $\gamma_z(s_e) = z\cdot s_e$ for all $e \in E^1$ and $v \in E^0$. The gauge-invariant ideal structure of $C^*(E)$ was studied in \cite{BHRS}. Here Bates \emph{et al.} identified the \emph{saturated, hereditary} subsets $H$ of $E^0$, and showed that a large class of gauge-invariant ideals in $C^*(E)$ correspond to subgraphs $r^{-1}(H) \subset E$ where $H$ is saturated and hereditary; the ideal associated to $r^{-1}(H)$ is denoted $I_H$, and contains the Cunz-Krieger algebra $C^*(r^{-1}(H))$ as a full corner. Ideally, the quotient $C^*(E) / I_H$ would be isomorphic to the Cuntz-Krieger algebra $C^*(s^{-1}(E^0 \setminus H))$ of the complementary subgraph. In fact, to realise $C^*(E)/I_H$ as a Cuntz-Krieger algebra, one needs to append sources (that is, vertices $v$ such that $r^{-1}(v)$ is empty) to $s^{-1}(E^0 \setminus H)$ to obtain what is referred to in \cite{BHRS} as the \emph{quotient graph} $E/H$. Using the uniqueness theorems for Cuntz-Krieger algebras, Bates \emph{et al.} show that $C^*(E) / I_H$ is canonically isomorphic to the Cuntz-Krieger algebra $C^*(E/H)$ \cite[Proposition~3.4]{BHRS}, and thereby identify the remainder of the gauge-invariant ideals in $C^*(E)$ \cite[Theorem~3.6]{BHRS}. They also produce a condition on $E$ under which all ideals of $C^*(E)$ are gauge-invariant \cite[Corollary~3.8]{BHRS}.

In recent work \cite[Section~3]{MT}, Muhly and Tomforde study the \emph{relative graph algebras} $C^*(E,V)$ of directed graphs $E$ using a construction which once again involved appending sources to $E$, and show that the Cuntz-Krieger algebra $C^*(E/H)$ of the quotient graph is canonically isomorphic to a relative graph algebra associated to $s^{-1}(E^0 \setminus H)$. 

For higher-rank graphs, the situation is more complicated. A higher-rank graph $\L$ can be thought of as a graph in which the paths have a shape or \emph{degree} in $\NN^k$ rather than a length in $\NN$. Associated to each higher-rank graph $\L$ there is a $C^*$-algebra $C^*(\L)$ generated by partial isometries associated to paths in $\L$ and carrying a strongly continuous gauge action $\gamma$ of $\TT^k$. The decompositions of a path in $\L$ must be in bijective correspondence with the decompositions of its degree in $\NN^k$; this is called the \emph{factorisation property}. The factorisation property poses significant complications for an analysis of the gauge-invariant ideals of $C^*(\L)$ using methods like those of \cite{BHRS}. The point is that it is not clear how to generalise the quotient graph construction from \cite{BHRS} to the higher-rank setting: because of the factorisation property, the addition of a source locally will have global effects on the higher-rank graph, so it is unclear how to reconcile multiple such operations.

In this paper we analyse the relative Cuntz-Krieger algebras associated to higher-rank graphs $\L$ with a view to studying the gauge-invariant ideal structure of $C^*(\L)$. Since the analysis of relative graph algebras in \cite{MT} involves appending sources to graphs, we would face the same difficulties in generalising it to the higher-rank setting as we would face in generalising the quotient graph construction of \cite{BHRS}. Instead, we study the relative Cuntz-Krieger algebras of finitely aligned higher-rank graphs by regarding them as universal objects generated by families of partial isometries. Our main objective is to establish versions of the gauge-invariant uniqueness theorem and the Cuntz-Krieger uniqueness theorem for relative Cuntz-Krieger algebras associated to higher-rank graphs, and we achieve these aims in Theorem~\ref{thm:GIUT} and Theorem~\ref{thm:CKUT}. The motivation for this is that the r\^oles played by $C^*(E/H)$ and the usual uniqueness theorems for graph algebras in \cite{BHRS} will be filled by a relative Cuntz-Krieger algebra associated to $s^{-1}(\L^0\setminus H)$ and the uniqueness theorems Theorem~\ref{thm:GIUT} and Theorem~\ref{thm:CKUT} in an analysis of the gauge-invariant ideal structure of $C^*(\L)$ for a finitely aligned $k$-graph $\L$.

In Section~\ref{sec:k-graphs and repns}, we given the definition of a $k$-graph and establish the notation we will need in later sections. In Section~\ref{sec:relCKalgs}, we associate a relative Cuntz-Krieger algebra $C^*(\L;\Ee)$ to each pair $\L,\Ee$ where $\L$ is a finitely aligned $k$-graph, and $\Ee$ is a collection of \emph{finite exhaustive} subsets of $\L$. We establish the existence of the \emph{core} subalgebra $C^*(\L;\Ee)^\gamma$ which is the fixed-point algebra for the gauge action, and adapt the methods of \cite[Section~3]{RSY2} to show that $C^*(\L;\Ee)^\gamma$ is AF. In Section~\ref{sec:thecore}, we say what it means for a collection $\Ee$ of finite exhaustive sets to be \emph{satiated}, and for such $\Ee$ we use the description of $C^*(\L;\Ee)^\gamma$ obtained in Section~\ref{sec:relCKalgs} to establish elementary conditions on a relative Cuntz-Krieger $(\L;\Ee)$-family under which it determines an injective homomorphism of $C^*(\L;\Ee)^\gamma$. In Section~\ref{sec:satiations}, we show how to produce from an arbitrary collection $\Ee$ of finite exhaustive sets an enveloping collection $\overline\Ee$ such that $\overline\Ee$ is satiated and $C^*(\L;\Ee) = C^*(\L;\overline\Ee)$. In Section~\ref{sec:uniqueness theorems}, we prove versions of the gauge-invariant and Cuntz-Krieger uniqueness theorems for $C^*(\L;\Ee)$ when $\Ee$ is satiated; the results of Section~\ref{sec:satiations} show how to apply these theorems to $C^*(\L;\Ee)$ when $\Ee$ is not satiated.

\section{Higher-rank graphs} \label{sec:k-graphs and repns}
We regard $\NN^k$ as an additive semigroup with identity 0. For $m,n \in \NN^k$, we write $m \vee n$ for their coordinate-wise maximum and $m \wedge n$ for their coordinate-wise minimum.

\begin{dfn} \label{dfn:k-graph} 
Let $k \in \NN \setminus \{0\}$.  A \emph{graph of rank $k$}, or \emph{$k$-graph}, is a pair $(\L,d)$ where $\L$ is a countable category and $d$ is a functor from $\L$ to $\NN^k$ which satisfies the \emph{factorisation property}:  For all $\lambda \in \Mor(\L)$ and all $m,n \in \NN^k$ such that $d(\lambda) = m+n$, there exist unique morphisms $\mu$ and $\nu$ in $\Mor(\L)$ such that $d(\mu) = m$, $d(\nu) = n$ and $\lambda = \mu\nu$.
\end{dfn}

Since we are regarding $k$-graphs as generalised graphs, we refer to elements of $\Mor(\L)$ as \emph{paths} and to elements of $\Obj(\L)$ as \emph{vertices} and we write $r$ and $s$ for the codomain and domain maps.

The factorisation property allows us to identify $\Obj(\L)$ with $\{\lambda \in \Mor(\L) : d(\lambda) = 0\}$. So we write $\lambda \in \L$ in place of $\lambda \in \Mor(\L)$, and when $d(\lambda) = 0$, we regard $\lambda$ as a vertex of $\L$.

Given $\lambda \in \L$ and $E \subset\L$, we define $\lambda E :=\{\lambda\mu : \mu \in E, r(\mu) = s(\lambda)\}$ and $E\lambda := \{\mu\lambda : \mu \in E, s(\mu) = r(\lambda)\}$.  In particular if $d(v) = 0$, then $v$ is a vertex of $\L$ and $vE = \{\lambda \in E : r(\lambda) = v\}$; similarly, $Ev = \{\lambda \in \L : s(\lambda) = v\}$.  We write
\[
\L^n := \{\lambda \in \L : d(\lambda) = n\}.
\]

The factorisation property ensures that if $l \le m \le n \in \NN^k$ and if $d(\lambda) = n$, then there exist unique paths denoted $\lambda(0, l)$, $\lambda(l, m)$ and $\lambda(m,n)$ such that $d(\lambda(0,l)) = l$, $d(\lambda(l,m)) = m - l$, and $d(\lambda(m,n)) = n-m$ and such that $\lambda = \lambda(0,l)\lambda(l,m)\lambda(m,n)$.

Given $k \in \NN\setminus\{0\}$, and $k$-graphs $(\L_1, d_1)$ and $(\L_2, d_2)$, we call a covariant functor $x : \L_1 \to \L_2$ a \emph{graph morphism} if it satisfies $d_2 \circ x = d_1$.

\begin{dfn} \label{dfn:common extensions} 
Let $(\L,d)$ be a $k$-graph. Given $\mu,\nu \in \L$, we say that $\lambda$ is a \emph{minimal common extension} of $\mu$ and $\nu$ if $d(\lambda) = d(\mu) \vee d(\nu)$, $\lambda(0, d(\mu)) = \mu$, and $\lambda(0, d(\nu)) = \nu$. We denote the collection of all minimal common extensions of $\mu$ and $\nu$ by $\MCE(\mu,\nu)$. We write $\Lmin(\mu,\nu)$ for the collection
\[
\Lmin(\mu,\nu) := \{(\alpha,\beta) \in \L \times \L : \mu\alpha =  \nu\beta \in \MCE(\mu,\nu)\}.
\]
If $E \subset \L$ and $\mu \in \L$, then we write $\Ext(\mu;E)$ for the set
\[
\Ext(\mu;E) := \{\alpha \in s(\mu)\L : (\alpha,\beta) \in \Lmin(\mu,\nu)\text{ for some }\nu \in E\}
\]
of extensions of $\mu$ with respect to $E$. We say that $\L$ is finitely aligned if $\MCE(\mu,\nu)$ is finite (possibly empty) for all $\mu,\nu \in \L$. 

Let $v \in \L^0$ and $E \subset v\L$. We say that $E$ is \emph{exhaustive} if $\Ext(\lambda;E)$ is nonempty for all $\lambda \in v\L$.
\end{dfn}

\begin{lemma} \label{lem:exhaustive sets} 
Let $(\L,d)$ be a finitely aligned $k$-graph, let $v \in \L^0$, let $E \subset v \L$ be finite and exhausitve, and let $\mu \in v \L$.  Then $\Ext(\mu; E)$ is a finite exhaustive subset of $s(\mu) \L$. Moreover $\mu \in E\L$ if and only if $s(\mu) \in \Ext(\mu;E)$.
\end{lemma}
\begin{proof} 
Let $E' := \Ext(\mu;E)$. Since $E$ is finite and $\L$ is finitely aligned we know that $E'$ is finite, and $E' \subset s(\mu)\L$ by definition, so we need only check that $E'$ is exhaustive.  Let $\sigma \in s(\mu)\L$.  Since $E$ is exhaustive, there exists $\lambda \in E$ with $\Lmin(\lambda,\mu\sigma) \not= \emptyset$, say $(\alpha,\beta) \in \Lmin(\lambda,\mu\sigma)$.  So $\lambda\alpha = \mu\sigma\beta$. Setting $\tau := (\mu\sigma\beta)(d(\mu), d(\lambda)\vee d(\mu))$, we have $\tau \in \Ext(\mu;\{\lambda\}) \subset E'$ by the factorisation property, and $\mu\sigma\beta = \mu\tau\tau'$ for some $\tau'$. But then the factorisation property gives $\sigma\beta = \tau\tau'$, so $(\sigma\beta)(0, d(\sigma)\vee d(\tau)) \in \MCE(\sigma,\tau)$. Since $\sigma \in s(\mu)\L$ was arbitrary, it follows that $E'$ is exhaustive. The last statement of the lemma follows from the factorisation property.
\end{proof}

\section{Relative Cuntz-Krieger algebras}\label{sec:relCKalgs}

\begin{notation} Let $(\L,d)$ be a finitely aligned $k$-graph. We define
\[\textstyle
\FE(\L) := \bigcup_{v \in \L^0}\{E \subset v\L\setminus\{v\} :\text{ $E$ is finite and exhaustive}\}.
\]
For $E \in \FE(\L)$ we write $r(E)$ for the vertex $v \in \L^0$ such that $E \subset v\L$.
\end{notation}

\begin{dfn}\label{dfn:relCK family}
Let $(\L,d)$ be a finitely aligned $k$-graph, and let $\Ee$ be a subset of $\FE(\L)$. A \emph{relative Cuntz-Krieger $(\L;\Ee)$-family} is a collection $\{t_\lambda : \lambda \in \L\}$ of partial isometries in a $C^*$-algebra satisfying
\begin{itemize}
\item[(TCK1)] $\{t_v : v \in \L^0\}$ is a collection of mutually orthogonal projections;
\item[(TCK2)] $t_\lambda t_\mu = t_{\lambda\mu}$ whenever $s(\lambda) = r(\mu)$;
\item[(TCK3)] $t^*_\lambda t_\mu = \sum_{(\alpha,\beta) \in \Lmin(\lambda,\mu)} t_\alpha t^*_\beta$ for all $\lambda,\mu \in \L$; and
\item[(CK)] $\prod_{\lambda \in E} (t_{r(E)} - t_\lambda t^*_\lambda) = 0$ for all $E \in \Ee$.
\end{itemize}
\end{dfn}

\begin{rmk}\label{rmk:pairwise commute}
Relation~(CK) is well-defined because (TCK3) ensures that the projections $\{t_\lambda t^*_\lambda : \lambda \in \L\}$ pairwise commute. Note also that (TCK3) together with the $C^*$-identity show that $t_v \not= 0$ for all $v \in \L^0$ if and only if $t_\lambda \not = 0$ for all $\lambda \in \L$.
\end{rmk}

For each finitely aligned $k$-graph $\L$, and each subset $\Ee$ of $\FE(\L)$ there exists a $C^*$-algebra $C^*(\L;\Ee)$ generated by a relative Cuntz-Krieger $(\L;\Ee)$-family $\{s_\Ee(\lambda) : \lambda \in \L\}$ which is universal in the sense that if $\{t_\lambda : \lambda \in \L\}$ is a relative Cuntz-Krieger $(\L;\Ee)$-family in a $C^*$-algebra $B$, then there exists a unique homomorphism $\pi^\Ee_t : C^*(\L;\Ee) \to B$ such that $\pi^\Ee_t(s_\Ee(\lambda)) = t_\lambda$ for all $\lambda \in \L$.

For $z = (z_1,\dots,z_k) \in \TT^k$ and $n = (n_1,\dots,n_k) \in \NN^k$, we write $z^n$ for the product $\prod^k_{i=1} z_i^{n_i} \in \TT$. The universal property of $C^*(\L;\Ee)$ guarantees that there exists a strongly continuous \emph{gauge action} $\gamma$ of $\TT^k$ on $C^*(\L;\Ee)$ such that $\gamma_z(s_\Ee(\lambda)) = z^{d(\lambda)} s_\Ee(\lambda)$ for all $\lambda \in \L$. Averaging over this gauge action gives a faithful conditional expectation $\Phi^\gamma_\Ee$ from $C^*(\L;\Ee)$ to the fixed point algebra 
\[
C^*(\L;\Ee)^\gamma = \clsp\{s_\Ee(\lambda) s_\Ee(\mu)^* : d(\lambda) = d(\mu)\}.
\]
We refer to $C^*(\L;\Ee)^\gamma$ as the \emph{core} of $C^*(\L;\Ee)$. The remainder of this section is devoted to showing that $C^*(\L;\Ee)^\gamma$ is AF. This material is adapted directly from \cite[Section~3]{RSY2}.

Recall from \cite{RSY2} that for a finitely aligned $k$-graph $(\L,d)$ and a finite subset $E \subset \L$, the set $\Pi E$ is the smallest subset of $\L$ such that $E \subset \Pi E$ and such that
\begin{equation}
\label{eq:prodclosed}
\lambda,\mu,\sigma \in G\text{ with } d(\lambda) = d(\mu)\text{ and } s(\lambda) = s(\mu) \text{ implies } 
\lambda\Ext(\mu;\{\sigma\}) \subset G.
\end{equation}
We write $\Pi E \times_{d,s} \Pi E$ for the set $\{(\lambda,\mu) \in \Pi E \times \Pi E : d(\lambda) = d(\mu), s(\lambda) = s(\mu)\}$. 

Lemma~3.2 of \cite{RSY2} shows that $\Pi E$ is finite, that if $\lambda, \mu \in \Pi E \times_{d,s} \Pi E$, then for $\nu \in s(\lambda) \L$, we have $\lambda\nu \in \Pi E$ if and only if $\mu\nu \in \Pi E$, and that if $\lambda,\mu \in \Pi E$ and $(\alpha,\beta) \in \Lmin(\lambda,\mu)$, then $\lambda\alpha \in \Pi E$.

\begin{dfn}\label{dfn:matrix algebras}
Let $(\L,d)$ be a finitely aligned $k$-graph, let $E \subset \L$ be finite, let $\Ee \subset \FE(\L)$ and let $\{t_\lambda :\lambda \in \L\}$ be a relative Cuntz-Krieger $(\L;\Ee)$-family. We define 
\[
M^t_{\Pi E} := \lsp\{t_\lambda t^*_\mu : (\lambda,\mu) \in \Pi E \times_{d,s} \Pi E\}.
\]
For $\lambda,\mu \in \Pi E \times_{d,s} \Pi E$, we define
\begin{equation}\label{eq:Theta def}
\Theta(t)^{\Pi E}_{\lambda,\mu} := 
t_\lambda \Big(\prod_{\lambda\nu \in \Pi E\setminus\{\lambda\}} (t_{s(\lambda)} - t_{\nu} t^*_{\nu})\Big) t^*_\mu.
\end{equation}
\end{dfn}

It is straightforward to check that Lemmas 3.2~and~3.11 and Proposition~3.9 of \cite{RSY2} apply to any family of partial isometries satisfying (TCK1)--(TCK3); for details see \cite[Chapter~3]{ASPhD}. Hence each $M^t_{\Pi E}$ is a finite-dimensional $C^*$-algebra \cite[Lemma~3.2]{RSY2}. Moreover, for all $(\lambda,\mu), (\sigma,\tau) \in \Pi E \times_{d,s} \Pi E$, we have 
\begin{equation}\label{eq:matrix units}
\big(\Theta(t)^{\Pi E}_{\lambda,\mu}\big)^* = \Theta(t)^{\Pi E}_{\mu,\lambda}
\quad\text{ and }\quad
\Theta(t)^{\Pi E}_{\lambda,\mu} \Theta(t)^{\Pi E}_{\sigma,\tau} = \delta_{\mu,\sigma} \Theta(t)^{\Pi E}_{\lambda,\tau}
\end{equation}
by \cite[Proposition~3.9]{RSY2}, and for $\lambda,\mu \in \Pi E \times_{d,s} \Pi E$, we have
\begin{equation}\label{eq:Thetas span}\textstyle
t_\lambda t^*_\mu = \sum_{\lambda\nu \in \Pi E} \Theta(t)^{\Pi E}_{\lambda\nu,\mu\nu}
\end{equation}
by \cite[Lemma~3.11]{RSY2}.

We can now show that the core is AF, and give a condition under which a representation of $C^*(\L;\Ee)$ is faithful on the core.

\begin{prop} \label{prp:faithful on core}
Let $(\L,d)$ be a finitely aligned $k$-graph, and let $\Ee$ be a subset of $\FE(\L)$. Then $C^*(\L;\Ee)^\gamma$ is an AF algebra. If $\{t_\lambda : \lambda \in \L\}$ is a relative Cuntz-Krieger $(\L;\Ee)$-family, then $\pi^\Ee_t$ is injective on $C^*(\L;\Ee)^\gamma$ if and only if $\Theta(t)^{\Pi E}_{\lambda,\mu}$ is nonzero whenever $\Theta(s_\Ee)^{\Pi E}_{\lambda,\mu}$ is nonzero.
\end{prop}
\begin{proof}
We have $C^*(\L;\Ee)^\gamma = \bigcup_{E \subset \L\text{ finite}} M^{s_\Ee}_{\Pi E}$, so $C^*(\L;\Ee)^\gamma$ is an AF algebra. If $E$ is a finite subset of $\L$, then~\eqref{eq:matrix units} shows that $\{\Theta(s_\Ee)^{\Pi E}_{\lambda,\mu} : (\lambda,\mu) \in \Pi E \times_{d,s} \Pi E, \Theta(s_\Ee)^{\Pi E}_{\lambda,\mu} \not= 0\}$ is a collection of nonzero matrix units, and~\eqref{eq:Thetas span} shows that these matrix units span $M^{s_\Ee}_{\Pi E}$. So if $\Theta(t)^{\Pi E}_{\lambda,\mu}$ is nonzero whenever $\Theta(s_\Ee)^{\Pi E}_{\lambda,\mu}$ is nonzero, then $\pi^\Ee_t$ is injective on each $M^{s_\Ee}_{\Pi E}$, and hence on $C^*(\L;\Ee)^\gamma$ by \cite[Lemma~1.3]{ALNR}. 
\end{proof}

\section{Nonzero matrix units and the $\Ee$-compatible boundary path representation} \label{sec:thecore}
In this section, we identify the \emph{satiated} subsets of $\FE(\L)$, and when $\Ee$ is satiated, we characterise the $\Theta(s_\Ee)^{\Pi E}_{\lambda,\mu}$ which are nonzero in $C^*(\L;\Ee)$.

\begin{dfn}\label{dfn:satiated}
Let $(\L,d)$ be a finitely aligned $k$-graph. We say that a subset $\Ee \subset \FE(\L)$ is \emph{satiated} if it satisfies
\begin{itemize}
\item[(S1)] If $G \in \Ee$ and $E \in \FE(\L)$ with $G \subset E$ then $E \in \Ee$.
\item[(S2)] If $G \in \Ee$ with $r(G) = v$ and if $\mu \in v\L \setminus G\L$, then $\Ext(\mu;G) \in \Ee$.
\item[(S3)] If $G \in \Ee$ and $0 < n_\lambda \le d(\lambda)$ for each $\lambda \in G$, then $\{\lambda(0, n_\lambda) : \lambda \in G\} \in \Ee$.
\item[(S4)] If $G \in \Ee$, $G' \subset G$, and $G'_\lambda \in \Ee$ with $r(G_\lambda') = s(\lambda)$ for each $\lambda \in G'$, then $\textstyle\big((G\setminus G') \cup \big(\bigcup_{\lambda \in G'} \lambda G'_\lambda\big)\big) \in \Ee$.
\end{itemize}
\end{dfn}

The remainder of this section is devoted to proving the following theorem.

\begin{theorem} \label{thm:characterise faithful reps}
Let $(\L,d)$ be a finitely aligned $k$-graph, and suppose that $\Ee \subset \FE(\L)$ is satiated. Let $\{t_\lambda : \lambda \in \L\}$ be a relative Cuntz-Krieger $(\L;\Ee)$-family. The homomorphism $\pi^\Ee_t$ is injective on $C^*(\L;\Ee)^\gamma$ if and only if
\begin{itemize}
\item[(1)] $t_v \not= 0$ for all $v \in \L^0$, and
\item[(2)] $\prod_{\lambda \in F}(t_{r(F)} - t_\lambda t^*_\lambda) \not= 0$ for all $F \in \FE(\L) \setminus \Ee$.
\end{itemize}
\end{theorem}

To prove Theorem~\ref{thm:characterise faithful reps}, we construct a representation of $C^*(\L;\Ee)$ which satisfies conditions (1)~and~(2). As usual, we obtain this representation by definining an appropriate boundary-path space. 

Recall from \cite{RSY1} that for $k \in \NN$ and $m \in (\NN \cup\{\infty\})^k$, the $k$-graph $\Omega_{k,m}$ has vertices $\{n \in \NN^k : n \le m\}$, morphisms $\{(n_1, n_2) : n_1, n_2 \in \NN^k, n_1 \le n_2 \le m\}$, degree map $d((n_1, n_2)) = n_2 - n_1$ and range and source maps $r((n_1, n_2)) = n_1$, $s((n_1, n_2)) = n_2$.

\begin{dfn}\label{dfn:rel bdry path} 
Let $(\L,d)$ be a finitely aligned $k$-graph, and let $\Ee \subset \FE(\L)$ be satiated. We say that a graph morphism $x: \Omega_{k, m} \to \L$ is an \emph{$\Ee$-compatible boundary path} of $\L$ if for every $n \in \NN^k$ such that $n \le m$, and every $E \in \Ee$ such that $r(E) = x(n)$, there exists $\lambda \in E$ such that $x(n, n+d(\lambda))  = \lambda$.  We denote the collection of all $\Ee$-compatible boundary paths of $\L$ by $\partial(\L;\Ee)$. We write $d(x)$ for $m$ and $r(x)$ for $x(0)$.
\end{dfn}

If $x \in \partial(\L;\Ee)$ and $\lambda \in \L r(x)$ then there is a unique graph morphism $\lambda x : \Omega_{k, d(\lambda) + d(x)} \to \L$ such that $(\lambda x)(0, d(\lambda)) = \lambda$ and $(\lambda x)(d(\lambda), n + d(\lambda)) = x(0, n)$ for all $n \le d(x)$. Likewise, if $n \in \NN^k$ with $n \le d(x)$, there is a unique graph morphism $x|^{d(x)}_n : \Omega_{k, d(x) - n} \to \L$ such that $x|^{d(x)}_n(0, m) = x(n, n+m)$ whenever $n+m \le d(x)$. These two constructions are inverse to each other in the sense that 
\begin{equation}\label{eq:extend/restrict}
(\lambda x)|^{d(\lambda x)}_{d(\lambda)} = x = \big(x(0,n)\big)\big(x|^{d(x)}_n\big)
\quad\text{for all $\lambda \in \L r(x)$ and all $n \le d(x)$.}
\end{equation}

\begin{lemma}\label{lem:bps ext/rest}
Let $(\L,d)$ be a finitely aligned $k$-graph, and let $\Ee \subset \FE(\L)$ be satiated. Let $x \in \partial(\L;\Ee)$. If $n \in \NN^k$ with $n \le d(x)$, then $x|^{d(x)}_n \in \partial(\L;\Ee)$, and if $\lambda \in \L r(x)$, then $\lambda x \in \partial(\L;\Ee)$. 
\end{lemma}
\begin{proof} 
For the first statement, just note that each vertex on $x|^{d(x)}_n$ is also a vertex on $x$. For the second statement, suppose $n\in \NN^k$ with $n \le d(\lambda x)$, and suppose $E \in \Ee$ with $r(E) = (\lambda x)(n)$.  Let $\lambda' = (\lambda x)(n, n \vee d(\lambda))$, and let $x' = x|_{(n\vee d(\lambda)) - d(\lambda)}^{d(x)}$, so that $(\lambda x)|_n^{d(\lambda x)} = \lambda' x'$, and $x' \in \partial(\L;\Ee)$ by the first statement of the lemma. We must show that there exists $\mu \in E$ such that $(\lambda' x')(0,d(\mu)) = \mu$.  If there exists $\mu \in E$ with $d(\mu) \le d(\lambda')$ and $\lambda'(0, d(\mu)) = \mu$, we are done, so we may assume that $\lambda' \not \in E\L$.  By (S2), we have $\Ext(\lambda',E) \in \Ee$, and $r(\Ext(\lambda',E)) = s(\lambda') = r(x')$ by definition.  Since $x' \in \partial(\L;\Ee)$, it follows that there exists $\alpha \in \Ext(\lambda';E)$ such that $x'(0,d(\alpha)) = \alpha$; equivalently, there exists $\mu \in E$ and $(\alpha,\beta) \in \Lmin(\lambda',\mu)$ such that $\alpha = x'(0,d(\alpha))$.  But now $\lambda'\alpha = \mu\beta$, and in particular, $(\lambda' x')(0,d(\mu)) = \big(\lambda' x'(0,d(\alpha))\big)(0,d(\mu)) = (\mu\beta)(0,d(\mu)) = \mu$.
\end{proof}

\begin{dfn}\label{dfn:rel bdry path repn} 
Let $(\L,d)$ be a finitely aligned $k$-graph and let $\Ee \subset \FE(\L)$ be satiated. Define partial isometries $\{S_{\Ee}(\lambda) : \lambda \in \L\} \subset \Bb(\ell^2(\partial(\L;\Ee)))$ by
\[
S_{\Ee}(\lambda) e_x := 
\begin{cases} e_{\lambda x} &\text{if $s(\lambda) = r(x)$} \\
0 &\text{otherwise.}
\end{cases}
\]
\end{dfn}

\begin{lemma}
Let $(\L,d)$ be a finitely aligned $k$-graph and let $\Ee \subset \FE(\L)$ be satiated. The collection $\{S_{\Ee}(\lambda) : \lambda \in \L\}$ is a relative Cuntz-Krieger $(\L;\Ee)$-family which we call the \emph{$\Ee$-compatible boundary-path representation} of $\L$.
\end{lemma}
\begin{proof}
First notice that Lemma~\ref{lem:bps ext/rest} ensures that for $\lambda \in \L$ and $x \in \partial(\L;\Ee)$, we have
\begin{equation}\label{eq:SsubEe adjoints}
S_\Ee(\lambda)^* e_x = 
\begin{cases}
e_{x|^{d(x)}_{d(\lambda)}} &\text{if $x(0,d(\lambda)) = \lambda$} \\
0 &\text{otherwise.}
\end{cases}
\end{equation}

For $v \in \L^0$, we have that $S_\Ee(v)$ is the projection onto $\clsp\{e_x : x \in v\partial(\L;\Ee)\}$ and hence $\{S_\Ee(v) : v \in \L^0\}$ are mutually orthogonal projections, establishing (TCK1). Since composition in the category $\L$ is associative, (TCK2) is straightforward to check. To see (TCK3) one uses \eqref{eq:SsubEe adjoints} to apply both $S_\Ee(\lambda)^* S_\Ee(\mu)$ and $\sum_{(\alpha,\beta) \in \Lmin(\lambda,\mu)} S_\Ee(\alpha) S_\Ee(\beta)^*$ to an arbitrary basis element $e_x$; calculations like those of \cite[Example~7.4]{RS1} show that the two agree. Finally, for condition (CK) let $x \in \partial(\L;\Ee)$ and $E \in \Ee$ with $r(E) = r(x) = x(0)$. Then $x(0,n) \in E$ for some $n \le d(x)$ by definition of $\partial(\L;\Ee)$, and we have $\big(S_\Ee(x(0)) - S_\Ee(x(0,n)) S_\Ee(x(0,n))^*\big) e_x = 0$ by \eqref{eq:extend/restrict}~and~\eqref{eq:SsubEe adjoints}. Since $\big(S_\Ee(x(0)) - S_\Ee(x(0,n)) S_\Ee(x(0,n))^*\big)$ is a term in $\prod_{\lambda \in E}\big(S_\Ee(r(E)) - s_\Ee(\lambda) s_\Ee(\lambda)^*\big)$, it follows that the kernel of the latter contains $e_x$. Since $x \in \partial(\L;\Ee)$ and $E \in \Ee$ were arbitrary, this establishes (CK).
\end{proof}

\begin{lemma} \label{lem:vDelta nonempty} 
Let $(\L,d)$ be a finitely aligned $k$-graph, let $\Ee$ be a satiated subset of\/ $\FE(\L)$, and let $v \in \L^{0}$. Then 
\begin{itemize}
\item[(1)] $v \partial(\L;\Ee)$ is nonempty. 
\item[(2)] If $F \in v\FE(\L) \setminus \Ee$, then $v\partial(\L;\Ee) \setminus F\partial(\L;\Ee)$ is nonempty.
\end{itemize}
\end{lemma}

To prove Lemma~\ref{lem:vDelta nonempty}, we first need the following technical lemma.

\begin{lemma}\label{lem:more in Ebar} 
Let $(\L,d)$ be a finitely aligned $k$-graph, and suppose that $\Ee \subset \FE(\L)$ is satiated. Suppose that $E \in \Ee$ and that $F \subset r(E)\L \setminus \L^0$ is finite and satisfies
\[
\mu \in E \implies 
\begin{cases}
\mu \in F\L &\text{or} \\
\Ext(\mu;F) \in \Ee.
\end{cases}
\]
Then $F \in \Ee$.
\end{lemma}
\begin{proof} Define $G := E\setminus F\L$, and for each $\mu \in G$, let $G_\mu := \Ext(\mu;F)$.  Then each $G_\mu \in \Ee$ by hypothesis, so (S4) gives
\[
\textstyle E' := \big((E \setminus G) \cup \big(\bigcup_{\mu \in G} \mu G_\mu\big)\big) \in \Ee.
\]
For $\lambda \in E \setminus G$ we have $\lambda(0,n) \in F$ for some n; in this case, let $n_\lambda := n$. Since $E \in \Ee$ we have $E \cap \L^0 = \emptyset$ and hence $n_\lambda > 0$. For $\mu \in G$ and $\lambda \in \mu G_\mu$, we have $\lambda = \mu\beta$ for some $\beta \in \Ext(\mu;F)$, so there exists $\sigma \in F$ and $\alpha \in \L$ such that $(\alpha,\beta) \in \Lmin(\sigma,\mu)$.  Hence
\[
\lambda(0, d(\sigma)) = (\mu\beta)(0, d(\sigma)) = (\sigma\alpha)(0,d(\sigma)) = \sigma \in F;
\]
in this case, set $n_\lambda : = d(\sigma)$. Since $F \cap \L^0 = 0$ by hypothesis, we have that $n_\lambda > 0$. Now $E'' := \{\lambda(0, n_\lambda) : \lambda \in E'\} \subset F$.  But $E' \in \Ee$, and hence (S3) ensures that $E'' \in \Ee$.  Since $E'' \subset F$ and $F$ is finite, it now follows from (S1) that $F \in \Ee$.
\end{proof}

To prove Lemma~\ref{lem:vDelta nonempty}, we also need the following a result due to Farthing, Muhly and Yeend.

\begin{lemma}[Farthing, Muhly and Yeend, 2003]\label{lem:FMY}
Let $(\L,d)$ be a $k$-graph. For $v \in \L^0$, $E \subset v\L$, $\lambda_1 \in v\L$, and $\lambda_2 \in s(\lambda_1)\L$, we have
$\Ext(\lambda_2; \Ext(\lambda_1; E)) = \Ext(\lambda_1\lambda_2; E)$.
\end{lemma} 
\begin{proof}
The result is proved in \cite{FMY}, currently in draft form; a proof also appears in \cite[Appendix~A]{ASPhD}.
\end{proof}

\begin{proof}[Proof of Lemma~\ref{lem:vDelta nonempty}]
The proofs of both statements of Lemma~\ref{lem:vDelta nonempty} proceed by constructing an $\Ee$-relative boundary path with the desired properties. The two constructions have a great deal in common, but the construction for statement~(2) is somewhat more complicated. To avoid duplication, we present the full text of the proof of statement~(2) below, but we typeset those parts of the proof which are germane only to statement~(2) in slanted text, and enclose them in square brackets \big[{\sl like this}\big]\!.

Define $P : (\NN \setminus \{0\})^2 \to (\NN \setminus \{0\})$ by
\[
P(m,n) := \frac{(m+n-1)(m+n-2)}{2} + m.
\]
Then $P$ is the position function corresponding to the diagonal listing of $(\NN\setminus\{0\})^2$ in the sense that if $(m,n)$ is the $l^{\rm th}$ term in the diagonal listing, then $P(m,n) = l$. For all $l \in \NN\setminus\{0\}$, define $(i_l, j_l) := P^{-1}(l)$. Fix $v \in \L^0$ {\big[}{\sl and fix $F \in v\FE(\L)\setminus\Ee$\/\,}{\big]}\!.

We claim that there exist a sequence $\{\lambda_l : l \ge 1\} \subset v\L$ and listings $\{E_{l,j} : j \ge 1\}$ of $s(\lambda_l)\Ee$ for all $l \ge 1$ satisfying 
\begin{itemize}
\item[(i)] $\lambda_{l+1}(0, d(\lambda_l)) = \lambda_l$ for all $l \ge 1$,
\item[(ii)] $\lambda_{l+1}\big(d(\lambda_{i_l}), d(\lambda_{l+1})\big)$ belongs to $E_{i_l, j_l} \L$ for all $l \ge 1$.
\item[{\big[\sl(iii)}] {\sl $\Ext(\lambda_{l+1}; F)$ belongs to $\FE(\L) \setminus\Ee$ for all $l \ge 0$.}\big]
\end{itemize}
We prove the claim by induction on $l$. For a basis case, let $l = 0$ and define $\lambda_{l+1} = \lambda_1 := v$. For each $w \in \L^0$, the collection of finite subsets of $w\L$ is countable because $\L$ is countable. In particular, $w \Ee$ is countable. Let $\{E_{1,j} : j \in \NN\setminus\{0\}\}$ be any listing of $v\Ee$. Note that (i)~and~(ii) are trivial in this case because $l = 0$ \big[{\sl and (iii) is satisfied because $\Ext(v;F) = F$}\big].

Now suppose as an inductive hypothesis that $l \ge 1$, and that $\lambda_n$ and $\{E_{n,j} : j \ge 1\}$ exist and satisfy (i)~and~(ii) \big[{\sl and (iii)}\big] for $1 \le n \le l$.

Let $\lambda^{l}_{i_l} := \lambda_l\big(d(\lambda_{i_l}), d(\lambda_{l})\big)$. Notice that $i_l < l$, so $E_{i_l, j_l} \in s(\lambda_{i_l})\Ee$ has already been defined by the inductive hypothesis. Suppose first that $\lambda^l_{i_l}$ belongs to $E_{i_l,j_l}\L$. Define $\lambda_{l+1} := \lambda_{l}$, and $E_{i_{l+1}, j} := E_{i_l, j}$ for all $j \ge 1$. We have that $\lambda_{l+1}$ satisfies~(i) by definition, and satisfies~(ii) because we supposed $\lambda^l_{i_l}$ to belong to $E_{i_l,j_l}\L$. \big[{\sl We have that $\lambda_{l+1}$ satisfies~(iii) because $\lambda_l$ satisfies~(iii) by the inductive hypothesis.}\big]

Now suppose that $\lambda^l_{i_l}$ does not belong to $E_{i_l,j_l}\L$. Let $E := \Ext(\lambda^l_{i_l}; E_{i_l, j_l})$. Then $E \not= \emptyset$ because $E_{i_l, j_l} \in \FE(\L)$. For $\alpha \in E$ we have that $\lambda^l_{i_l}\alpha = \mu\beta$ for some $\mu \in E$ and $(\alpha,\beta) \in \Lmin(\lambda^l_{i_l}, \mu)$. It follows that for any $\nu_{l+1}$ in $E$ we have that $\lambda_{l+1} := \lambda_{l}\nu_{l+1}$ satisfies~(ii). Such a choice of $\lambda_{l+1}$ trivially satisfies~(i).

\ \hskip-1em\big[{\sl To complete the construction of $\lambda_{l+1}$, we need only show that there exists a choice of $\nu_{l+1} \in E$ such that $\lambda_{l+1} := \lambda_{l} \nu_{l+1}$ also satisfies~(iii). Since $\lambda_{l}$ satisfies (iii), we have that $F_{l} := \Ext(\lambda_{l}; F)$ belongs to $\FE(\L) \setminus \Ee$. By the contrapositive of Lemma~\ref{lem:more in Ebar}, there exists $\alpha \in E\setminus F_{l}\L$ such that $\Ext(\alpha;F_{l}) \not\in \Ee$. But Lemma~\ref{lem:exhaustive sets} ensures that $\Ext(\alpha;F_{l}) \in \FE(\L)$, so $\Ext(\alpha; F_{l}) \in \FE(\L) \setminus \Ee$. Let $\nu_{l+1} := \alpha$, and define $\lambda_{l+1} := \lambda_{l}\nu_{l+1}$. Then 
\begin{equation}\label{eq:satisfies (ii)}
\Ext(\lambda_{l+1}; F) = \Ext(\lambda_{l}\nu_{l+1}; F) = \Ext(\nu_{l+1}; \Ext(\lambda_{l}; F))
\end{equation}
by Lemma~\ref{lem:FMY}. But $\Ext(\lambda_{l}; F) = F_{l}$ by definition, so \eqref{eq:satisfies (ii)} gives $\Ext(\lambda_{l+1}; F) = \Ext(\nu_{l+1}; F_{l})$ which belongs to $\FE(\L)\setminus\Ee$ by choice of $\nu_{l+1}$. Hence $\lambda_{l+1}$ satisfies~(iii) as required.\big]}

Let $m := \lim_{l \to \infty} d(\lambda_l) \in (\NN \cup \{\infty\})^k$. Since $\{\lambda_l : l \ge 1\}$ satisfies~(i), there exists a unique graph morphism $x : \Omega_{k, m} \to \L$ such that $x(0, d(\lambda_l)) = \lambda_{l}$ for all $l \in \NN\setminus\{0\}$. 

We have that $r(x) = v$ by definition, so to see that $x \in v\partial(\L;\Ee)$, suppose that $M \in \NN^k$ with $M \le m$. Let $E \in x(M)\Ee$. We must show that there exists $N \ge M$ such that $x(M,N) \in E$. By definition of $x$ there exists $l \ge 1$ such that $M \le d(\lambda_l)$. If $\lambda_{l}(M, d(\lambda_{l}))$ belongs to $E\L$, then we are done, so suppose that $\lambda_{l}(M, d(\lambda_{l})) \not\in E\L$. By (S3), it follows that $G := \Ext(\lambda_{l}(M, d(\lambda_{l})); E) \in s(\lambda_{l})\Ee$, and hence that $G = E_{i_l, j}$ for some $j \ge 1$. But then property~(ii) ensures that $\lambda_{P(i_l,j) + 1}(M,N) \in E$ for some $N$, and it follows that $x(M,N) \in E$ as required.

\ \hskip-1em\big[{\sl Finally we must show that $x \not\in F\L$. Suppose for contradiction that $x \in F\L$. Then $x(0,N) \in F$ for some $N$, and it follows from the definition of $x$ that there exists $l \ge 1$ such that $\lambda_{l}(0,N) = x(0,N) \in F$. Hence $s(\lambda_{l})$ belongs to $\Ext(\lambda_{l}; F)$. But for $G \in \FE(\L)$, we have $G \cap \L^0 = \emptyset$ by definition, so $s(\lambda_{l}) \in \Ext(\lambda_l; F)$ contradicts~(iii). Hence $x \not\in F\L$.}\big]
\end{proof}

\begin{cor}\label{cor:which MUs nonzero}
Let $(\L,d)$ be a finitely aligned $k$-graph, and let $\Ee \subset \FE(\L)$ be satiated. The vertex projections $s_\Ee(v)$ are all nonzero. Moreover, if $E \subset v\L\setminus\L^0$ is finite, then $\prod_{\lambda \in E}(s_\Ee(v) - s_\Ee(\lambda) s_\Ee(\lambda)^*) = 0$ if and only if $E \in \Ee$. 
\end{cor}

To prove Corollary~\ref{cor:which MUs nonzero}, we make use of an equality established in \cite{RSY2}: let $(\L,d)$ be a finitely aligned $k$-graph, let $\{t_\lambda : \lambda \in \L\}$ be a collection of partial isometries satisfying (TCK1)--(TCK3), let $v$ be an element of $\L^0$, let $E$ be a finite subset of $v\L$, and let $\mu$ be an element of $v\L$. Then \cite[Equation~(3.4)]{RSY2} shows that
\begin{equation}\label{eq:shift gaps}
\Big(\prod_{\lambda \in E} (t_v - t_\lambda t^*_\lambda)\Big) t_\mu t^*_\mu 
= t_\mu\Big(\prod_{\alpha \in \Ext(\mu;E)} (t_{s(\mu)} - t_\alpha t^*_\alpha)\Big)t^*_\mu,
\end{equation} 
with the convention that the empty product is equal to the unit of the multiplier algebra so that if $\Ext(\mu;E) = \emptyset$ then the right-hand side of~\eqref{eq:shift gaps} is equal to $t_\mu t^*_\mu$.

\begin{proof}[Proof of Corollary~\ref{cor:which MUs nonzero}]
Statement (1) of Lemma~\ref{lem:vDelta nonempty} shows that for $v \in \L^0$, there exists $x \in v\partial(\L;\Ee)$, and then $S_\Ee(v)e_x = e_x \not= 0$. So $S_{\Ee}(v)$ is nonzero, and the universal property of $C^*(\L;\Ee)$ then shows that $s_{\Ee}(v)$ is nonzero.

For the second statement of the Corollary, the ``if" direction is precisely~(CK). For the reverse implication, suppose that $E \subset v\L$ but $E \not\in \Ee$. If $E \not\in \FE(\L)$ then there exists $\xi \in v\L$ such that $\Ext(\xi;E) = \emptyset$. Equation~\ref{eq:shift gaps} shows that
\[
s_{\Ee}(\xi) s_{\Ee}(\xi)^* \prod_{\lambda \in E} (s_\Ee(v) - s_\Ee(\lambda) s_\Ee(\lambda)^*) = s_{\Ee}(\xi) s_{\Ee}(\xi)^*,
\]
and hence $\prod_{\lambda \in E} (s_\Ee(v) - s_\Ee(\lambda) s_\Ee(\lambda)^*)$ is nonzero by Remark~\ref{rmk:pairwise commute}. On the other hand, if $E \in \FE(\L)\setminus\Ee$, then statement~(2) of Lemma~\ref{lem:vDelta nonempty} shows that there exists $x \in v\partial(\L;\Ee)\setminus E\partial(\L;\Ee)$. We then have $\prod_{\lambda \in E} (S_\Ee(v) - S_\Ee(\lambda) S_\Ee(\lambda)^*) e_{x} = e_{x} \not= 0$, and then the universal property of $C^*(\L;\Ee)$ gives $\prod_{\lambda \in E} (s_\Ee(v) - s_\Ee(\lambda) s_\Ee(\lambda)^*) \not= 0$.
\end{proof}

\begin{proof}[Proof of Theorem~\ref{thm:characterise faithful reps}]
The ``only if" implication follows from Corollary~\ref{cor:which MUs nonzero}. Equation~\eqref{eq:Theta def} and Corollary~\ref{cor:which MUs nonzero} show that for $E \subset \L$ finite and $(\lambda,\mu) \in \Pi E \times_{d,s} \Pi E$, we have $\Theta(s_\Ee)^{\Pi E}_{\lambda,\mu} = 0$ if and only if $T^{\Pi E}(\lambda) := \{\nu \in \L \setminus\L^0 : \lambda\nu \in \Pi E\}$ belongs to $\Ee$. Hence for the ``if" direction it suffices to establish that if $\Theta(t)^{\Pi E}_{\lambda,\mu} = 0$, then $T^{\Pi E}(\lambda)$ belongs to $\Ee$; indeed, by~\eqref{eq:Theta def}, it suffices to show that if $T^{\Pi E}(\lambda) \not\in \Ee$, then $\prod_{\nu \in T^{\Pi E}(\lambda)}(t_{s(\lambda)} - t_\nu t^*_\nu) \not= 0$. 

So suppose $T^{\Pi E}(\lambda) \not \in \Ee$. If $T^{\Pi E}(\lambda) \not\in \FE(\L)$, then $\prod_{\nu \in T^{\Pi E}(\lambda)}(t_{s(\lambda)} - t_\nu t^*_\nu) \not= 0$ exactly as in the proof of Corollary~\ref{cor:which MUs nonzero}. On the other hand, if $T^{\Pi E}(\lambda) \in \FE(\L) \setminus \Ee$, then $\prod_{\nu \in T^{\Pi E}(\lambda)}(t_{s(\lambda)} - t_\nu t^*_\nu) \not= 0$ by assumption.
\end{proof}

\section{Constructing satiations}\label{sec:satiations}
In this section we show how to use Theorem~\ref{thm:characterise faithful reps} to characterise the homomorphisms of arbitrary relative Cuntz-Krieger algebras $C^*(\L;\Ee)$ which are injective on the core, and not just those for which $\Ee$ is satiated.

\begin{dfn}
Let $(\L,d)$ be a finitely aligned $k$-graph, and let $\Ee \subset \FE(\L)$. We write $\overline\Ee$ for the smallest satiated subset of $\FE(\L)$ which contains $\Ee$, and we call $\overline\Ee$ the \emph{satiation} of $\Ee$.
\end{dfn}

The idea is to show that for any $\Ee \subset \FE(\L)$, we have $C^*(\L;\Ee) = C^*(\L;\overline\Ee)$. To this end we define maps $\Sigma_1$--$\Sigma_4$ on subsets of $\FE(\L)$, and show that iterated application of these maps produces $\overline\Ee$ from $\Ee$.

\begin{dfn}
Let $(\L,d)$ be a finitely aligned $k$-graph, and for $\Ee \subset \FE(\L)$, define
\begin{align*}
\Sigma_1(\Ee) &= \{F \subset \L\setminus\L^0 : F\text{ is finite, and there exists } E \in \Ee \text{ with } E \subset F\} \\
\Sigma_2(\Ee) &= \{\Ext(\mu;E) : E \in \Ee, \mu \in r(E)\L \setminus E\L\} \\
\Sigma_3(\Ee) &= \big\{\{\lambda(0, n_\lambda) : \lambda \in E\} : E \in \Ee, 0 < n_\lambda \le d(\lambda) \text{ for all } \lambda \in E\big\} \\
\Sigma_4(\Ee) &\textstyle = \big\{(E \setminus F) \cup\big(\bigcup_{\lambda \in F} \lambda F_\lambda\big) : E \in \Ee, F \subset E, \\ 
&\hskip5cm F_\lambda \in s(\lambda)\Ee\text{ for all } \lambda \in F\big\}.
\end{align*}
\end{dfn}

\begin{lemma}\label{lem:SigmaSubIs}
Let $(\L,d)$ be a finitely aligned $k$-graph, and let $\Ee \subset \FE(\L)$. Then $\Ee \subset \Sigma_i(\Ee) \subset \FE(\L)$ for $1 \le i \le 4$. Let $\{t_\lambda : \lambda \in \L\}$ be a relative Cuntz-Krieger $(\L;\Ee)$-family and let $E \in \Sigma_i(\Ee)$ for $1 \le i \le 4$. Then $\prod_{\lambda \in E}(t_{r(E)} - t_\lambda t^*_\lambda) = 0$.
\end{lemma} 
\begin{proof}
Let $E \in \Ee$. We trivially have $E \in \Sigma_1(\Ee)$. To see that $E \in \Sigma_2(\Ee)$, note that $r(E) \not\in E\L$ by definition, and $E = \Ext(r(E); E)$. To see that $E \in \Sigma_3(\Ee)$, just take $n_\lambda := d(\lambda)$ for all $\lambda \in E$. Finally, to see that $E \in \Sigma_4(\Ee)$, take $F = \emptyset\subset E$.

We will now establish that if $\{t_\lambda : \lambda \in \L\}$ is a relative Cuntz-Krieger $(\L;\Ee)$-family and $E \in \Sigma_i(\Ee)$, then $\prod_{\lambda \in E}(t_{r(E)} - t_\lambda t^*_\lambda) = 0$. 

If $i = 1$, then  $E = G \cup F$ for some $G \in \Ee$ and finite $F \subset r(G)\L$, and $\prod_{\lambda \in E} (t_{r(E)} - t_\lambda t^*_\lambda) \le \prod_{\lambda \in G} (t_{r(G)} - t_\lambda t^*_\lambda) = 0$.

If $i = 2$, then $E = \Ext(\mu; G)$ for some $G\in\Ee$ and $\mu \in r(G)\L\setminus G\L$. So multiplying~\eqref{eq:shift gaps} by $t_\mu^*$ on the left and by $t_\mu$ on the right gives
\[
\prod_{\lambda \in E} (t_{r(E)} - t_\lambda t^*_\lambda) 
= t^*_\mu \Big(\prod_{\sigma \in G} (t_{r(G)} - t_\sigma t^*_\sigma)\Big) t_\mu = 0. 
\]

If $i = 3$, then $E = \{\lambda(0, n_\lambda) : \lambda \in G\}$ for some $G \in \Ee$ and $0 < n_\lambda \le d(\lambda)$  for each $\lambda \in G$. Since $t_{r(E)} - t_{\lambda(0,n_\lambda)} t^*_{\lambda(0, n_\lambda)} \le t_{r(E)} - t_\lambda t^*_\lambda$ for all $\lambda \in E$, we then have
\[
\prod_{\lambda \in E} (t_{r(E)} - t_\lambda t^*_\lambda) \le \prod_{\mu \in G} (t_{r(G)} - t_\mu t^*_\mu) = 0.
\]

If $i = 4$, then $E = G \setminus G' \cup\big(\bigcup_{\lambda \in G'} \lambda G'_\lambda\big)$ for some $G \in \Ee$, $G' \subset G$, and $G'_\lambda \in s(\lambda)\Ee$ for each $\lambda \in G'$. Lemma~C.7 of \cite{RSY2} shows that for $\lambda \in G'$, we have $t_{r(G)} - t_\lambda t^*_\lambda = \prod_{\mu \in G'_\lambda} (t_{r(G)} - t_{\lambda\mu} t^*_{\lambda\mu})$.  Hence
\begin{align*} 
\prod_{\lambda \in E} (t_{r(E)} - t_\lambda t^*_\lambda) 
&= \Big(\prod_{\lambda \in G \setminus G'} (t_{r(E)} - t_\lambda t^*_\lambda)\Big) \prod_{\lambda \in G'} \Big(\prod_{\mu \in G'_\lambda} (t_{r(E)} - t_{\lambda\mu} t^*_{\lambda\mu}) \Big) \\
&= \prod_{\lambda \in G} (t_{r(G)} - t_\lambda t^*_\lambda).
\end{align*}

It remains only to show that $\Sigma_i(\Ee) \subset \FE(\L)$. For this, first notice that $E \in \Sigma_i(\Ee)$ implies that $E \cap \L^0 = \emptyset$ and that $E$ is finite by definition of $\Sigma_1$--$\Sigma_4$. Now let $\{t_\lambda : \lambda \in \L\}$ be a relative Cuntz-Krieger $(\L;\Ee)$-family in which $t_v \not= 0$ for all $v \in \L^0$; such a family exists by Corollary~\ref{cor:which MUs nonzero}. Suppose that $v \in \L^0$ and that $E$ is a finite subset of $v\L\setminus\L^0$ with $\prod_{\lambda \in E}(t_v - t_\lambda t^*_\lambda) = 0$, and suppose for contradiction that $E \not\in \FE(\L)$. Then there exists $\mu \in v\L$ such that $\Lmin(\mu;\lambda) = \emptyset$ for all $\lambda \in E$. Equation~\ref{eq:shift gaps} gives $t_\mu t^*_\mu \prod_{\lambda \in E}(t_v - t_\lambda t^*_\lambda) = t_\mu t^*_\mu$, and hence $t_\mu t^*_\mu = 0$, contradicting $t_v \not = 0$ for all $v \in \L^0$. Since we have already established that if $\{t_\lambda : \lambda \in \L\}$ is a relative Cuntz-Krieger $(\L;\Ee)$-family and $E \in \Sigma_i(\Ee)$, then $\prod_{\lambda \in E}(t_{r(E)} - t_\lambda t^*_\lambda) = 0$, it follows that $\Sigma_i(\Ee) \subset \FE(\L)$ as required.
\end{proof}

\begin{notation}
We write $\Sigma$ for the map $\Sigma_4 \circ \Sigma_3 \circ \Sigma_2 \circ \Sigma_1$. For $n \in \NN$ and $\Ee \subset \FE(\L)$, we write $\Sigma^n(\Ee)$ for
\[
\overbrace{\Sigma\circ\Sigma\circ\cdots\circ\Sigma}^{n\text{ terms}}(\Ee),
\]
and write $\Sigma^{\infty}(\Ee)$ for $\bigcup^\infty_{n=1} \Sigma^n(\Ee)$.
\end{notation}

\begin{prop}\label{prp:Ebar=SigmaE}
Let $(\L,d)$ be a finitely aligned $k$-graph and let $\Ee \subset \FE(\L)$. Then $\Sigma^\infty(\Ee) = \overline\Ee$.
\end{prop}
\begin{proof}
The definitions of the maps $\Sigma_1$--$\Sigma_4$ show that $\Sigma_i(\Ee) \subset \overline\Ee$ for all $i$, and hence that $\Sigma^\infty(\Ee) \subset \overline\Ee$. Hence, it suffices to show that $\Sigma^\infty(\Ee)$ is satiated. If $G \in \Sigma^\infty(\Ee)$ and $E$ is constructed from $G$ as in (S1), (S2) or (S3), then we have $G \in \Sigma^n(\Ee)$ for some $n \in \NN$, and then since we have $\Ee \subset \Sigma_i(\Ee)$ for all $i$ by Lemma~\ref{lem:SigmaSubIs}, it follows that $E \in \Sigma^{n+1}(\Ee) \subset \Sigma^\infty(\Ee)$ as required. If $F \in \Sigma^\infty(\Ee)$, $G \subset F$, and $G_\lambda \in s(\lambda)\Sigma^\infty(\Ee)$ for all $\lambda \in G$, then there exist $n \in \NN$ with $F \in \Sigma^n(\Ee)$, and $n_\lambda \in \NN$ such that $G_\lambda \in \Sigma^{n_\lambda}(\Ee)$ for each $\lambda \in \L$. Let $N := \max\{n,n_\lambda : \lambda \in G\}$. Again since Lemma~\ref{lem:SigmaSubIs} shows that $\Ee \subset \Sigma_i(\Ee)$ for all $i$, we have that $F$ and each $G_\lambda$ belong to $\Sigma^N(\Ee)$. The definition of $\Sigma_4$ together with another application of Lemma~\ref{lem:SigmaSubIs} shows that $E \in \Sigma^{N+1}(\Ee) \subset \Sigma^{\infty}(\Ee)$, and the proof is complete.
\end{proof}

\begin{cor}
Let $(\L,d)$ be a finitely aligned $k$-graph, and let $\Ee$ be any subset of $\FE(\L)$. Then $C^*(\L;\Ee) = C^*(\L;\overline\Ee$).
\end{cor}
\begin{proof}
An induction on $n$ using the last statement of Lemma~\ref{lem:SigmaSubIs} shows that if $F \in \Sigma^n(\Ee)$, then $\prod_{\mu \in F}(t_{r(F)} - t_\mu t^*_\mu) = 0$ for all $n \in \NN$. Hence $\prod_{\mu \in F}(t_{r(F)} - t_\mu t^*_\mu) = 0$ for all $F \in \overline\Ee$ by Proposition~\ref{prp:Ebar=SigmaE}. It follows that every relative Cuntz-Krieger $(\L;\Ee)$-family is a relative Cuntz-Krieger $(\L;\overline\Ee)$-family. On the other hand $\Ee \subset \overline\Ee$ by definition, so every Cuntz-Krieger $(\L;\overline\Ee)$-family is trivially a relative Cuntz-Krieger $(\L;\Ee)$-family. The universal properties of $C^*(\L;\Ee)$ and $C^*(\L;\overline\Ee)$ now show that the two algebras coincide.
\end{proof}

\section{Uniqueness theorems}\label{sec:uniqueness theorems}
In this section we prove versions of the gauge-invariant uniqueness theorem and the Cuntz-Krieger uniqueness theorem for $C^*(\L;\Ee)$.

\begin{theorem}\label{thm:GIUT}
Let $(\L,d)$ be a finitely aligned $k$-graph, and let $\Ee \subset \FE(\L)$ be satiated. Let $\{t_\lambda : \lambda \in \L\}$ be a relative Cuntz-Krieger $(\L;\Ee)$-family in a $C^*$-algebra $B$, and suppose that
\begin{itemize}
\item[(1)] $t_v \not= 0$ for all $v \in \L^0$;
\item[(2)] $\prod_{\lambda \in F}(t_{r(F)} - t_\lambda t^*_\lambda) \not= 0$ for all $F \in \FE(\L) \setminus \Ee$; and
\item[(3)] there exists an action $\theta : \TT^k \to \Aut(B)$ such that $\theta_z(t_\lambda) = z^{d(\lambda)} t_\lambda$ for all $z \in \TT^k$ and $\lambda \in \L$.
\end{itemize}
Then $\pi^\Ee_t$ is injective.
\end{theorem}
\begin{proof}
Theorem~\ref{thm:characterise faithful reps} and Conditions (1)~and~(2) guarantee that $\pi^\Ee_t$ is injective on $C^*(\L;\Ee)^\gamma$. Assume without loss of generality that $B = C^*(\{t_\lambda : \lambda \in \L\})$. Since the polynomials are continuous on $\TT^k$, and since $B = \clsp\{t_\lambda t^*_\mu : \lambda,\mu \in \L\}$ by~(TCK3), we have that $\theta$ is strongly continuous. Since $\pi^\Ee_t$ is equivariant in $\theta$ and $\gamma$, averaging over $\theta$ gives a norm-decreasing linear map $\Phi^\theta_\Ee$ on $B$ which satisfies $\Phi^\theta_\Ee \circ \pi^\Ee_t = \pi^\Ee_t \circ \Phi^\gamma_\Ee$. The result now follows from an argument identical to that of \cite[Proposition~4.1]{RSY2}.
\end{proof}

To state our Cuntz-Krieger uniqueness theorem, we first need to establish some notation.

\begin{dfn}\label{dfn:MCE for BPs}
Let $(\L,d)$ be a $k$-graph, and let $x : \Omega_{k,d(x)} \to \L$ and $y : \Omega_{k,d(y)} \to \L$ be graph morphisms. We say that a graph morphism $z : \Omega_{k,d(z)} \to \L$ is a \emph{minimal common extension} of $x$ and $y$ if it satisfies
\begin{itemize}
\item[(1)] $d(z)_j = \max\{d(x)_j, d(y)_j\}$ for $1 \le j \le k$; and
\item[(2)] $z|_{\Omega_{k,d(x)}} = x$ and $z|_{\Omega_{k,d(y)}} = y$.
\end{itemize}
We write $\MCE(x,y)$ for the collection of minimal common extensions of $x$ and $y$.
\end{dfn}

It turns out that to obtain a Cuntz-Krieger uniqueness theorem for relative Cuntz-Krieger algebras, the appropriate analogue of an aperiodic path is a path $x \in \partial(\L;\Ee)$ such that
\begin{equation}\label{eq:x aperiodic}
\text{for distinct } \lambda,\mu \in \L r(x), \text{ we have } \MCE(\lambda x, \mu x) = \emptyset.
\end{equation}

\begin{theorem}\label{thm:CKUT}
Let $(\L,d)$ be a finitely aligned $k$-graph and let $\Ee \subset \FE(\L)$ be satiated. Suppose that $(\L,\Ee)$ satisfies
\begin{equation}\label{eq:new aperiodicity}
\begin{split}
&\text{For all $v \in \L^0$ there exists $x \in v\partial(\L;\Ee)$ satisfying \eqref{eq:x aperiodic},} \\
&\quad\text{and for all $v\in \L^0$ and $F \in v\FE(\L) \setminus \Ee$ there exists} \\
&\quad\text{$x \in v\partial(\L;\Ee)\setminus F\partial(\L;\Ee)$ satisfying \eqref{eq:x aperiodic}.}
\end{split} \tag{C}
\end{equation}
Let $\{t_\lambda : \lambda \in \L\}$ be a relative Cuntz-Krieger $(\L;\Ee)$-family such that $t_v \not= 0$ for all $v \in \L^0$, and $\prod_{\lambda \in F} (t_{r(F)} - t_\lambda t^*_\lambda) \not= 0$ for all $F \in \FE(\L)\setminus\Ee$. Then $\pi^\Ee_t$ is injective.
\end{theorem}

The remainder of the section is devoted to proving Theorem~\ref{thm:CKUT}. We first need some technical lemmas.

\begin{lemma}\label{lem:finite no MCE}
Let $(\L,d)$ be a finitely aligned $k$-graph, and suppose that $x : \Omega_{k,d(x)} \to \L$ is a graph morphism satisfying \eqref{eq:x aperiodic}. Suppose that $\lambda \not= \mu$ with $s(\lambda) = s(\mu) = r(x)$. Then there exists $n^x_{\lambda,\mu} \in \NN^k$ such that 
\[
n^x_{\lambda,\mu} \le d(x)
\quad\text{ and }\quad 
\Lmin(\lambda x(0, n^x_{\lambda,\mu}), \mu x(0, n^x_{\lambda,\mu})) = \emptyset.
\]
\end{lemma}
\begin{proof}
Suppose for contradiction that for all $n \in \NN^k$ with $n \le d(x)$, we have $\Lmin(\lambda x(0,n), \mu x(0,n)) \not= \emptyset$.

For each $i \in \NN$, define $n(i) \in \NN^k$ by $n(i)_j := \min\{d(x)_j, i\}$. By assumption, there exists $(\alpha_i, \beta_i) \in \Lmin(\lambda x(0, n(i)), \mu x(0, n(i)))$ for each $i \in \NN$; since $\Lmin(\lambda, \mu)$ is finite, there must exist a pair $(\eta_1, \zeta_1)$ belonging to $\Lmin(\lambda, \mu)$ and an infinite subset $I_1 \subset \NN$ such that for all $i \in I_1$,
\[
(\lambda x(0, n(i))\alpha_i)(0, d(\lambda) \vee d(\mu)) = \lambda\eta_1 = \mu\zeta_1.
\]
Set $i_1 := \min I_1$. For each $j \in I_1$ with $j > i_1$, we have
\[\begin{split}
\big(\lambda x(0, n(j))\alpha_j\big)&\big(0, d(\lambda x(0, n(i_1))) \vee d(\mu x(0, n(i_1)))\big) \\
&= \big(\lambda x(0, n(j))\alpha_j\big)\big(0, (d(\lambda) \vee d(\mu)) + n(i_1)\big) \\
&\in \MCE(\lambda x(0, n(i_1)), \mu x(0, n(i_1))).
\end{split}\]
Since $\MCE(\lambda x(0, n(i_1)), \mu x(0, n(i_1)))$ is finite, there exists a pair $(\eta_2, \zeta_2)$ belonging to $\Lmin(\lambda x(0, n(i_1)), \mu x(0, n(i_1)))$ and an infinite subset $I_2 \subset I_1\setminus\{i_1\}$ such that for each $i \in I_2$, we have
\[
\big(\lambda x(0, n(j))\alpha_j\big)\big(0, (d(\lambda) \vee d(\mu)) + n(i)\big) 
= \lambda x(0, n(i_1)) \eta_2 
= \mu x(0, n(i_1)) \zeta_2.
\]
Since $I_2 \subset I_1$, a straightforward calculation using the fact that $\lambda x(0, n(i_1)) \eta_2$ is an initial segment of $\lambda x(0, n(i)) \alpha_i$ for any $i \in I_2$ shows that
\[
(\lambda x(0, n(i_1)) \eta_2)(0, d(\lambda) \vee d(\mu)) = \lambda\eta_1 = \mu\zeta_1.
\]
Set $i_2 := \min I_2$. Iterating this procedure, we obtain a sequence 
\[
\{\sigma_l := \lambda x(0, n(i_l)) \eta_{l+1} : l \in \NN\}
\]
such that $d(\sigma_l) = (d(\lambda) \vee d(\mu)) + n(i_l)$, and $\sigma_{l+1}(0, d(\sigma_l)) = \sigma_l$ for all $l$. There is a unique graph morphism $y : \Omega_{k, d(y)} \to \L$ such that 
\[
d(y) = \lim_{l \to \infty} (d(\lambda) \vee d(\mu)) + n(i_l) = d(\lambda x) \vee d(\mu x),
\] 
and $y(0, d(\sigma_l)) = \sigma_l$ for all $l$. We then have
\[\begin{split}
y(0, d(\lambda) + n(i_l)) 
&= \sigma_l(0, d(\lambda) + n(i_l)) \\
&=\big(\lambda x(0, n(i_l)) \eta_{l+1}\big)(0, d(\lambda) + n(i_l))
= \lambda x(0, n(i_l)).
\end{split}\]
Since $n(i_l) \to d(x)$, it follows that $y|_{\Omega_{k, d(\lambda) + d(x)}} = \lambda x$. Similarly, $y|_{\Omega_{k, d(\mu) + d(x)}} = \mu x$. It follows that $y \in \MCE(\lambda x, \mu x)$, contradicting \eqref{eq:x aperiodic}.
\end{proof}

\begin{lemma}\label{lem:shift not in Ebar}
Let $(\L,d)$ be a finitely aligned $k$-graph, let $\Ee \subset \FE(\L)$ be satiated, and suppose that $F \in \FE(\L) \setminus \Ee$. Let $x \in r(F)\partial(\L;\Ee) \setminus F\partial(\L;\Ee)$, and let $n \in \NN^k$ with $n \le d(x)$. Then $\Ext(x(0,n); F) \in \FE(\L) \setminus \Ee$. 
\end{lemma}
\begin{proof}
By Lemma~\ref{lem:exhaustive sets}, we have $\Ext(x(0,n);F) \in \FE(\L)$. Suppose for contradiction that $\Ext(x(0,n);F) \in \Ee$. Since $x \in \partial(\L;\Ee)$, there exists $m > n$ such that $m \le d(x)$ and $x(n,m) \in \Ext(x(0,n);F)$. So there exists $\lambda \in F$ and $\alpha \in s(\lambda)\L$ such that $(\alpha, x(n,m)) \in \Lmin(\lambda, x(0,n))$. But then $x(0,m) = x(0,n)x(n,m) = \lambda\alpha$, contradicting the assumption that $x$ does not belong to $F\partial(\L;\Ee)$.
\end{proof}

\begin{cor}\label{cor:moving nonzero gaps}
Let $(\L,d)$ be a finitely aligned $k$-graph, let $\Ee \subset \FE(\L)$ be satiated, and suppose that $F \in \FE(\L) \setminus \Ee$. Let $x \in r(F)\partial(\L;\Ee) \setminus F\partial(\L;\Ee)$, and let $n \in \NN^k$ with $n \le d(x)$. Let $\{t_\lambda : \lambda \in \L\}$ be a relative Cuntz-Krieger $(\L;\Ee)$-family such that $t_v \not= 0$ for all $v \in \L^0$, and $\prod_{\lambda \in F} (t_{r(F)} - t_\lambda t^*_\lambda) \not= 0$ for all $F \in \FE(\L)\setminus\Ee$. Then 
\[
\prod_{\lambda \in F} (t_{r(F)} - t_\lambda t^*_\lambda) t_{x(0,n)} t^*_{x(0,n)}  = t_{x(0,n)} \Big(\prod_{\beta \in \Ext(x(0,n);F)} (t_{x(n)} - t_\beta t^*_\beta)\Big) t^*_{x(0,n)},
\]
and in particular is nonzero.
\end{cor}
\begin{proof}
The displayed equation is an instance of \eqref{eq:shift gaps}. Lemma~\ref{lem:shift not in Ebar} ensures that $\Ext(x(0,n);F) \in \FE(\L)\setminus \Ee$, and then $\prod_{\beta \in \Ext(x(0,n);F)} (t_{x(n)} - t_\beta t^*_\beta) \not= 0$ by hypothesis.
\end{proof}

\begin{lemma}\label{lem:tail gives FCE}
Let $(\L,d)$ be a finitely aligned $k$-graph, let $\Ee \subset \FE(\L)$ be satiated, and suppose that $(\L,\Ee)$ satisfies condition~\eqref{eq:new aperiodicity}. Let $\{t_\lambda : \lambda \in \L\}$ be a relative Cuntz-Krieger $(\L;\Ee)$-family such that $t_v \not= 0$ for all $v \in \L^0$, and $\prod_{\lambda \in F} (t_{r(F)} - t_\lambda t^*_\lambda) \not= 0$ for all $F \in \FE(\L)\setminus\Ee$. Let $\pi^\Ee_t$ be the representation of $C^*(\L;\Ee)$ determined by $\pi^\Ee_t(s_\Ee(\lambda)) = t_\lambda$. Let $a \in \lsp\{s_\Ee(\lambda) s_\Ee(\mu)^* : \lambda,\mu \in \L\} \subset C^*(\L;\Ee)$. Then $\|\pi^\Ee_t(\Phi^\gamma(a))\| \le \|\pi^\Ee_t(a)\|$.
\end{lemma}
\begin{proof}
Express $a = \sum_{\lambda,\mu \in \Pi E} a_{\lambda,\mu} s_\Ee(\lambda) s_\Ee(\mu)^*$ for some finite $E \subset \L$, and express $\Phi^\gamma(a) = \sum_{(\lambda,\mu) \in \Pi E \times_{d,s} \Pi E} b_{\lambda,\mu} \Theta(s_\Ee)^{\Pi E}_{\lambda,\mu}$; so we have
\[
\pi^\Ee_t(a) = \sum_{\lambda,\mu \in \Pi E} a_{\lambda,\mu} t_\lambda t_\mu^* \quad\text{and}\quad
\pi^\Ee_t(\Phi^\gamma(a)) = \sum_{(\lambda,\mu) \in \Pi E \times_{d,s} \Pi E} b_{\lambda,\mu} \Theta(t)^{\Pi E}_{\lambda,\mu}.
\]
Since the $\Theta(t)^{\Pi E}_{\lambda,\mu}$ are matrix units, there exists $n$ in $d(\Pi E)$ and $v \in s(\Pi E \cap \L^n)$ such that
\[
\|\pi^\Ee_t(\Phi^\gamma(a))\| = \Big\|\sum_{\lambda,\mu \in (\Pi E)v \cap \L^n} b_{\lambda,\mu} \Theta(t)^{\Pi E}_{\lambda,\mu}\Big\|.
\]
Write $T^{\Pi E}(n,v)$ for $T^{\Pi E}(\lambda)$ where $\lambda \in (\Pi E)v \cap \L^n$, so 
\[
T^{\Pi E}(n,v) = \{\nu \in \L\setminus\L^0 : \lambda\nu \in \Pi E\text{ for any }\lambda \in (\Pi E)v \cap \L^n\}.
\] 
Equation~\eqref{eq:prodclosed} ensures that $T^{\Pi E}_{n,v}$ is well-defined. If $T^{\Pi E}(n,v)$ belongs to $\Ee$, then we must have $\Phi^\gamma(a) = 0$ in which case the result is trivial. So suppose that $T^{\Pi E}(n,v) \not\in \Ee$. We claim that there exists $x \in v\partial(\L;\Ee) \setminus T^{\Pi E}(n,v)\partial(\L;\Ee)$ satisfying \eqref{eq:x aperiodic}. To see this, note that if $T^{\Pi E}(n,v) \in \FE(\L)$, then such an $x$ exists because $(\L,\Ee)$ satisfies condition~\eqref{eq:new aperiodicity}, whereas if $T^{\Pi E}(n,v) \not\in \FE(\L)$, then there exists $\sigma \in v\L$ with $\Ext(\sigma;T^{\Pi E}(n,v)) = \emptyset$, and condition~\eqref{eq:new aperiodicity} gives $x' \in s(\sigma)\L$ satisfying \eqref{eq:x aperiodic}; it is then easy to check that $x := \sigma x'$ also satisfies \eqref{eq:x aperiodic} and does not have an initial segment in $T^{\Pi E}(n,v)$.

For all $\lambda \in \Pi E$ with $d(\lambda) \le n$, $\mu \in (\Pi E) s(\lambda)$ with $\mu \not= \lambda$, and $\nu \in s(\lambda)\L$ such that $\lambda\nu \in (\Pi E)v \cap \L^n$, the factorisation property ensures that  $\lambda\nu \not= \mu\nu$. Hence Lemma~\ref{lem:finite no MCE} shows that there exists $n^x_{\lambda\nu,\mu\nu} \in \NN^k$ with $n^x_{\lambda\nu,\mu\nu} \le d(x)$ such that $\Lmin\big(\lambda\nu x(0,n^x_{\lambda\nu,\mu\nu}), \mu\nu x(0, n^x_{\lambda\nu,\mu\nu})\big) = \emptyset$. Define 
\[
N := \bigvee\{n^x_{\lambda\nu,\mu\nu} : \lambda,\mu \in \Pi E, d(\lambda) \not=d(\mu), \lambda\nu \in (\Pi E)v \cap \L^n\}.
\]
Since each $n^x_{\lambda\nu,\mu\nu} \le d(x)$, we have $N \le d(x)$, and for each $\lambda,\mu,\nu$ as above, we have that $\Lmin(\lambda\nu x(0,N), \mu\nu x(0,N)) = \emptyset$.

Define projections $P_1$ and $P_2$ by 
\[
P_1 := \sum_{\lambda \in (\Pi E)v \cap \L^n} \Theta(t)^{\Pi E}_{\lambda,\lambda} \quad\text{ and }\quad
P_2 := \sum_{\lambda \in (\Pi E)v \cap \L^n} t_{\lambda x(0,N)} t^*_{\lambda x(0,N)}.
\]
We have $P_1 \pi^\Ee_t(\Phi^\gamma(a)) = \sum_{\lambda,\mu \in (\Pi E)v \cap \L^n} b_{\lambda,\mu} \Theta(t)^{\Pi E}_{\lambda,\mu}$ and hence $\|P_1 \pi^\Ee_t(\Phi^\gamma(a))\| = \|\pi^\Ee_t(\Phi^\gamma(a))\|$. For $\lambda, \sigma \in (\Pi E)v \cap \L^n$, we have  $t_{\sigma x(0,N)}^* t_\lambda = \delta_{\sigma,\lambda} t^*_{x(0,N)}$ by (TCK3), and so for $\lambda,\mu \in (\Pi E)v \cap \L^n$, we have
\begin{align*}
P_2 \Theta(t)^{\Pi E}_{\lambda,\mu} P_2 
&= P_2 t_\lambda \Big(\prod_{\nu \in T^{\Pi E}(n,v)}(t_v - t_\nu t^*_\nu)\Big) t^*_\mu P_2 \\
&= t_{\lambda x(0,N)} t^*_{x(0,N)} \Big(\prod_{\nu \in T^{\Pi E}(n,v)}(t_v - t_\nu t^*_\nu)\Big) t_{x(0,N)}  t^*_{\mu x(0,N)} \\
&= t_{\lambda x(0,N)} \Big(\prod_{\beta \in \Ext(x(0,N); T^{\Pi E}(n,v))} t_{x(N)} - t_\beta t^*_\beta\Big) t^*_{\mu x(0,N)}
\end{align*}
by Corollary~\ref{cor:moving nonzero gaps}; Corollary~\ref{cor:moving nonzero gaps} also shows that this last expression is nonzero. For $\lambda \in (\Pi E)v \cap \L^n$, we have $t_{\lambda x(0,N)} \in \L^{n + N}$, and it follows that for $\lambda,\mu \in (\Pi E)v \cap \L^n$, we have $t^*_{\lambda x(0,N)} t_{\mu x(0,N)} =
\delta_{\lambda,\mu} t_{x(n)}$. Hence 
\[\textstyle
\big\{P_2 \Theta(t)^{\Pi E}_{\lambda,\mu} P_2 : \lambda,\mu \in (\Pi E)v \cap \L^n\big\}
\]
is a collection of nonzero matrix units, and compression by $P_2$ therefore implements an isomorphism of $M^t_{\Pi E}(n,v)$. It follows that 
\[
\|P_2 (P_1 \pi^\Ee_t(\Phi^\gamma(a))) P_2\| = \|P_1 \pi^\Ee_t(\Phi^\gamma(a))\| = \|\pi^\Ee_t(\Phi^\gamma(a))\|.
\]

On the other hand, we have $\|P_2 (P_1 \pi^\Ee_t(a)) P_2)\| \le \|\pi^\Ee_t(a)\|$ because $P_1$ and $P_2$ are projections. Thus, the proof of Lemma~\ref{lem:tail gives FCE} will be complete if we can establish that $P_2 (P_1 \pi^\Ee_t(a)) P_2) = P_2 (P_1 \pi^\Ee_t(\Phi^\gamma(a))) P_2)$. To do this, it suffices to show that if $\lambda,\mu \in \Pi E$ with $d(\lambda) \not= d(\mu)$ and $s(\lambda) = s(\mu)$, we have $P_2(P_1 t_\lambda t^*_\mu)P_2 = 0$. To see this, fix $\lambda,\mu \in \Pi E$ with $d(\lambda) \not= d(\mu)$ and $s(\lambda) = s(\mu)$, and calculate
\begin{align*}
P_2 P_1 t_\lambda t^*_\mu P_2
&= P_2\Big(\sum_{\lambda\nu \in (\Pi E)v \cap \L^n} 
      t_{\lambda\nu} \Big(\prod_{\substack{\lambda\nu\sigma' \in \Pi E \\ d(\sigma') > 0}} 
                   (t_{s(\nu)} - t_{\sigma'} t^*_{\sigma'})\Big) t^*_{\mu\nu}\Big) P_2 \\
&= \sum_{\lambda\nu \in (\Pi E)v \cap \L^n}\Big(
      \Big(\prod_{\substack{\lambda\nu\sigma' \in \Pi E \\ d(\sigma') > 0}} 
                   (t_{\lambda\nu}t^*_{\lambda\nu} - t_{\lambda\nu\sigma'} t^*_{\lambda\nu\sigma'})\Big)
                   P_2 t_{\lambda\nu} t^*_{\mu\nu} P_2\Big),
\end{align*}
because (TCK3) ensures that the projections $\{t_\lambda t^*_\lambda : \lambda \in \L\}$ pairwise commute. So it suffices to show that $P_2 t_{\lambda\nu}t^*_{\mu\nu} P_2 = 0$ for all $\nu$ such that $\lambda\nu \in (\Pi E)v \cap \L^n$. Fix such a $\nu$. We have that $\sigma,\tau \in \L^n$ implies $t^*_\sigma t_\tau = \delta_{\sigma,\tau} t_{s(\sigma)}$ by (TCK3). It follows that for $\sigma \in (\Pi E)v \cap \L^n$, we have 
\[
t^*_{\sigma x(0,N)} t_{\lambda\nu} = t^*_{x(0,N)} t^*_\sigma t_{\lambda\nu} = \delta_{\sigma,\lambda\nu} t^*_{x(0,N)}.
\]
Consequently, $P_2 t_{\lambda\nu} = t_{\lambda\nu x(0,N)} t^*_{x(0,N)}$. Hence
\begin{align*}
P_2 t_{\lambda\nu} t^*_{\mu\nu} P_2 
&= t_{\lambda\nu x(0,N)} t^*_{x(0,N)} t^*_{\mu\nu} P_2 \\
&= t_{\lambda\nu x(0,N)} \sum_{\tau \in (\Pi E)v \cap \L^n} t^*_{\mu\nu x(0,N)} t_{\tau x(0,N)} t^*_{\tau x(0,N)}.
\end{align*} 
Since $d(\mu) \not= d(\lambda)$, we have $d(\mu\nu) \not= d(\tau\nu)$, and hence $\mu\nu \not= \tau\nu$ for each $\tau \in (\Pi E)v \cap \L^n$. It follows that $\Lmin\big(\mu\nu x(0,N), \tau\nu x(0,N)\big) = \emptyset$ for all $\tau \in (\Pi E)v \cap \L^n$ by our choices of $x$ and $N$. Hence the final line of the above calculation is equal to zero by (TCK3), proving the Lemma.
\end{proof}

\begin{proof}[Proof of Theorem~\ref{thm:CKUT}]
Lemma~\ref{lem:tail gives FCE} shows that the formula 
\[
t_\lambda t^*_\mu \mapsto \delta_{d(\lambda), d(\mu)} t_\lambda t^*_\mu
\]
extends to a norm-decreasing linear map $\Phi^t$ on $\pi^\Ee_t(C^*(\L;\Ee))$. Replacing $\Phi^\theta$ with $\Phi^t$ in the proof Theorem~\ref{thm:GIUT} now establishes the result.
\end{proof}

\end{document}